\documentclass[a4paper, 11pt, margin=1in]{article}
\usepackage{amsmath,amsthm,amssymb,mathrsfs,graphicx,url,bbm}
\usepackage{fullpage}
\usepackage{multirow}
\usepackage[usenames]{color}
\usepackage{ifthen}
\usepackage{setspace}
\usepackage{caption}
\usepackage{subfigure}
\usepackage{float}
\usepackage{graphicx}
\usepackage{algorithm}
\usepackage{bm}
\usepackage[noend]{algpseudocode}
\usepackage{enumitem}
\usepackage{amsfonts}

\newtheorem{theorem}{Theorem}
\newtheorem{corollary}{Corollary}

\usepackage{natbib}
\usepackage{authblk}

\makeatletter
\def\BState{\State\hskip-\ALG@thistlm}
\makeatother

\def\keywords{\vspace{.5em}
{\textbf{Keywords}:\,\relax%
}}

\theoremstyle{plain}

\theoremstyle{remark}
\newtheorem{remark}{Remark}

\theoremstyle{definition}

\theoremstyle{remark}

\theoremstyle{definition}

\newcommand{\prob}{\mathsf{P}}

\renewcommand{\phi}{\varphi}

\title{Fisher information matrix of binary time series}
\author{Xu Gao$^{1}$, 
	Hernando Ombao$^{1,2,3}$\, 
	Daniel Gillen$^{1}$ \\
	$^{1}$Department of Statistics, University of California, Irvine, California, U.S.A. \\
	$^{2}$Department of Cognitive Sciences, University of California, Irvine, California, U.S.A.\\
	$^{3}$Program on Applied Mathematics \& Computational Science, \\
	King Abdullah University of Science and Technology, Saudi Arabia\\
\tt{xgao2@uci.edu \ \ hernando.ombao@kaust.edu.sa \ \ dgillen@uci.edu}
}
\date{}

\begin{document}

\maketitle
\begin{abstract}
A common approach to analyzing categorical correlated time series data is 
to fit a generalized linear model (GLM) with past data as covariate inputs. There 
remain challenges to conducting inference for short time series length. By treating 
the historical data as covariate inputs, standard errors of estimates of GLM 
parameters computed using the empirical Fisher information do not fully account 
the auto-correlation in the data. To overcome this serious limitation,
we derive the exact conditional Fisher information matrix of a general logistic autoregressive model with endogenous covariates for any series length $T$.  
Moreover, we also develop an iterative computational formula that allows for 
relatively easy implementation of the proposed estimator. Our simulation studies 
show that confidence intervals derived using the exact Fisher information matrix 
tend to be narrower than those utilizing the empirical Fisher information matrix 
while maintaining type I error rates at or below nominal levels.
Further, we establish that the exact Fisher information matrix approaches, as 
$T$ tends to infinity, the asymptotic Fisher information matrix previously derived 
for binary time series data. The developed exact conditional Fisher information matrix 
is applied to time-series data on respiratory rate among a cohort of expectant 
mothers where it is found to provide narrower confidence intervals for functionals 
of scientific interest and lead to greater statistical power when compared to the empirical Fisher information matrix.

\end{abstract}

\newpage
\section{Introduction}
\label{s:intro}

Time series data are widely collected in many fields such as
genetics, medicine and transportation(see \citet{Gouveia:2017}, \citet{Chen:2018A}, \citet{Chen:2018B}, \citet{Guo:2018}).
Various models for categorical time series that take into account temporal correlation are discussed in \citet{Kedem:1994},
\citet{Kedem:1980}, \citet{Diggle Liang:1994}, \citet{Fahrmeir:1994} and \citet{Gao:2017}, among others.
In this paper, we consider the logistic autoregressive model for binary time series 
data. Under this model, we derive the {\it exact} conditional Fisher information (Ex-FI) matrix for binary time series with arbitrary length $T$ and demonstrate that a correctly
specified Ex-FI leads to more efficient inference for regression parameters. 
In particular,  confidence intervals are narrower compared to those obtained 
using the empirical Fisher information (Em-FI) matrix \citep{Dodge:2003} 
while maintaining type I error rates at, or below, nominal levels. 

We briefly describe some of the related approaches to modeling binary time 
series. \citet{Keenan:1982} developed a model with an underlying unobserved 
process that is Gaussian first-order autoregressive. For binary time series with 
a Markovian structure, \citet{Billingsley:1961}, \citet{Meyn:2012}, 
\citet{Bonney:1987}, \citet{Fahrmeir:1987},  \citet{Kaufmann:1987}, 
\citet{Keenan:1982} and \ \citet{Muenz:1985} developed an inferential 
procedure based on the conditional likelihood. A comprehensive modeling 
framework based on partial likelihood inference and generalized linear models 
was developed in \citet{Fokianos:2003} and \citet{Kedem:2002}.

In practice, standard software for fitting generalized linear models (GLMs) to
binary time series 
use the past series values as ``explanatory variables" in the conditional mean 
of the response for the regression model \citep{Sybolt:1998}. One limitation 
of this approach is that it does not differentiate between explanatory variables 
that are exogenous to the time series data versus those that are endogenous 
(i.e., explanatory values that are past values of the time series). Thus, it does 
not properly take into account the auto-correlation structure in the data, leading 
to potentially undesirable consequences. In particular, the standard errors of 
the regression parameter estimates derived using the Em-FI matrix also ignore 
the auto-correlation structure. We demonstrate that this may lead to incorrect 
inference because the asymptotic covariance matrix of the partial 
(conditional) maximum likelihood estimators of the logistic model parameters are 
also incorrect. To exemplify the practical importance of this result, 
\citet{Sybolt:1998} utilized a logistic autoregressive model (LAR/LARX) to predict 
the outcome of supervised exercise for intermittent claudication. The inference 
did not  distinguish between covariates that were exogenous versus endogenous 
to the time series, hence yielding potentially invalid and/or inefficient 
statistical inference. Hence, our goal here is to derive the Ex-FI matrix
(for finite time series length $T$) of a logistic autoregressive (LAR/LARX) model 
for statistical inference. This model takes into account the correlation in binary 
time series data.


Asymptotic inference for modeling independent and serially-correlated binary 
responses (or binary time series) has been well studied in the literature. In the 
case of independent binary responses, several papers have established and 
discussed the poor performance of Wald-based tests and confidence intervals 
for the probability of success when utilizing the empirical information 
matrix (see \citet{Hauck:1977}; \citet{Newcombe:1998}; \citet{Agresti:2001}). 
As a result, the use of either the likelihood ratio or score based confidence intervals 
is generally recommended for inference in this setting. For binary time 
series,  \citet{Fahrmeir:1987} and \citet{Kaufmann:1987} used the 
Markovian assumption to demonstrate the asymptotic normality and efficiency 
of the maximum likelihood estimator under standard regularity 
conditions. \citet{Fokianos:1998} extended the idea by introducing time dependent
covariates (including past series values).  Under the framework of partial
likelihood inference, \citet{Fokianos:1998} proved the existence of an 
asymptotic conditional Fisher information (AFI) matrix and established 
asymptotic results of the maximum partial likelihood estimator. However, 
our simulation studies indicate that  one needs be cautious in applying 
asymptotic results when the length of the time series is small or moderate 
(say, $T < 200$). \citet{Startz:2012} provided a statistical strategy for modeling the binary autoregressive moving average (BARMA) under mild assumptions. 
In their study, one of the main considerations was the lack of analytical forms 
of the autocorrelation and the unconditional mean because of the nonlinearity 
of the model. In this paper, we propose a rigorous approach to derive the 
Ex-FI matrix of a LAR/LARX model that provides more efficient asymptotic 
inference in terms of narrow interval estimation while maintaining nominal 
type I error rate. 



In \cite{Fokianos:1998}, the AFI matrix was derived for the general case where 
the conditional distribution of a time series depends on its
own historical data as well as other covariates. While impressive in its generality, 
the primary limitation of this result is that it does not provide a closed form of 
the Fisher information matrix for specific models. As we will demonstrate, the form 
of the Fisher information matrix is non-trivial even for the logistic first order autoregressive LAR(1) model, perhaps the simplest model for binary time series. 
The difficulty in deriving the analytic form of the Fisher information matrix lies in the fact that the score function, or the Hessian matrix, contains cross-covariance related to time-varying covariates. Another limitation is that the Ex-FI matrix has not been derived for finite $T$. Instead,
only an asymptotic approximation based on the partial likelihood, which
turned out to be equivalent to the Em-FI matrix for the
LAR model, was provided. There are major consequences of these limitations. First, the result lacks the precise form of
the Fisher information matrix to conduct inference on specific LAR coefficients and functionals of these coefficients (e.g., probability of $Y_t = 1$ given the past values of
the binary series). Second, when $T$ is not sufficiently large, the discrepancy
between the Ex-FI and Em-FI matrices could lead to poor power, incorrect
significance level of tests, inefficient inference, and potentially misleading results from data analysis. {Third, the large sample theory derived in \citet{Kedem:2002} is based on the crucial assumption that $
	\frac{1}{T}\sum\limits_{t=1}^T\mathbb{I}(X_t\in A)\rightarrow\nu(A),
	$
	where $\nu(.)$ is a probability measure, $A$ is a Borel set and $\mathbb{I}(.)$ is the indicator function. Even when $T$ is large, such assumption may not be easily met. In this way, using Em-FI rather than Ex-FI may be misleading since no large sample theory is guaranteed.} 

Motivated by these limitations, this paper provides a derivation of the Ex-FI matrix of a
LAR/LARX model for arbitrary finite $T$. While the derivation is non-trivial we provide a computationally tractable expression
that can be easily implemented in an iterative manner. We report findings from simulation
studies suggesting that the derived Ex-FI matrix yields superior results relative to the Em-FI for small to moderate sample sizes. When compared to using the
Em-FI, inference based on the Ex-FI matrix produces narrower confidence intervals for a fixed significance
level; close to expected false positive rate and higher power when conducting tests of hypotheses. The simulation studies also demonstrate that the Ex-FI matrix converges
to the general AFI  developed in \citet{Fokianos:1998} in the sense that the norm of the
difference between the entries of the two matrices converges to 0 when the length $T$ of the binary time series
increases. Finally, we apply the developed Ex-FI matrix to time-series data on respiratory rate among a cohort of expectant mothers. Results show the similar pattern observed from simulations.  Namely, the Ex-FI matrix is found to provide narrower confidence intervals for functionals of scientific interest (such as the probability or log odds) and produce more statistical power when compared to the Em-FI matrix.

The remainder of this paper is organized as follows. In Section~\ref{section:Ex-FI}, we first derive the Ex-FI matrix of LAR/LARX model in general. We also propose a computation framework through functional iteration to obtain the Ex-FI matrix explicitly. At the end, we consider a special case when the order of LAR mode is 1 and calculate the analytic form of the Ex-FI matrix. In Section~\ref{section:sim}, we present some simulation results to compare the Ex-FI with Em-FI. Results show the benefit of using Ex-FI in terms of shorter confidence interval length and reasonable Type I error rate. Moreover, asymptotic behavior is also studied. In Section~\ref{section:real}, we applied the Ex-FI matrix to time-series data on respiratory rate among expectant mothers. By comparing with the Em-FI, we conclude that using Ex-FI can produce greater power and shorter confidence intervals when conducting statistical inference.

\section{Derivation of the Exact Conditional Fisher Information Matrix}
\label{section:Ex-FI}

\subsection{Logistic autoregressive model of order $p$ (LAR(p))}
\begin{theorem}
\label{thm1}
Consider a binary-valued correlated time series data $Y_t, \ t=1, \ldots, T$ where the conditional distribution of $Y_t$ depends on the
previous values via the conditional probability
\begin{eqnarray}\label{Eq:Pt}
P_t = \prob (Y_t=1\mid y_{t-1}, y_{t-2}, \cdots, y_1)=\frac{\exp (\bm{y_{-t}'\beta})}{1+\exp (\bm{y_{-t}'\beta})},
\end{eqnarray}
where $\bm{y_{-t}}=(1,y_{t-1},\cdots, y_{t-p})'$ is endogenous to the series and $\bm{\beta}=(\beta_0, \cdots, \beta_p)'.$ The exact Fisher Information (Ex-FI) matrix takes the form
\begin{equation}
I(\bm{\beta}\mid y_{p},\cdots, y_1)= \sum_{t=p+1}^{T}\sum\limits_{\bm{(y_{-t})}}\Big{[}\frac{\exp (\bm{y_{-t}'\beta})}{\{1+\exp (\bm{y_{-t}'\beta})\}^2}\bm{y_{-t}}\bm{y_{-t}}'\Big{]} Q_t (y_{t-1}, \cdots, y_{t-p}),
\end{equation}
where the conditional joint probability of $Y_{t-1}, \ldots, Y_{t-p}$ is
derived to be
\begin{eqnarray}
&&Q_t (y_{t-1}, \cdots, y_{t-p}) \nonumber\\ & = & \prob (Y_{t-1}=y_{t-1}, Y_{t-2}=y_{t-2}, \cdots, Y_{t-p}=y_{t-p}\mid y_{p}, y_{p-1}, \cdots, y_1) \nonumber \\
{} &= & \left\{
\begin{array}{lr}
\qquad 1, & \text{if}\quad  t=p+1 \\
\prod\limits_{k=1}^{t-p-1} P_{t-k}^{y_{t-k}}(1-P_{t-k})^{1-y_{t-k}}, & \text{if} \quad p+2\leq t \leq 2p+1 \\
\sum\limits_{(y_{p+1}, \cdots, y_{t-p-1})}\prod\limits_{k=1}^{t-p-1} P_{t-k}^{y_{t-k}}(1-P_{t-k})^{1-y_{t-k}}, & \text{if}\quad  t \geq 2p+2.
\end{array}
\right.
\label{marginal}
\end{eqnarray}
\end{theorem}

\textbf{Proof} The proof directly follows by the fact that the conditional log-likelihood function of $\bm{\beta}$ and the vector of conditional \
score functions are, respectively,
\begin{eqnarray*}
	\ell(\bm{\beta}\mid \bm{Y}) & = & \sum_{t=p+1}^{T}[Y_t(\bm{Y_{-t}'\beta})-\log\{1+\exp(\bm{Y_{-t}'\beta})\}] \\
	\bm{U(\beta, Y)} & = & \sum_{t=p+1}^{T}[\bm{Y_{-t}}\{Y_t-\frac{\exp (\bm{Y_{-t}'\beta})}{1+\exp (\bm{Y_{-t}'\beta})}\}].
\end{eqnarray*}
And
\begin{equation}
\frac{\partial ^2}{\partial \bm{\beta} \partial \bm{\beta^T}}\ell(\bm{\beta}\mid \bm{Y})= -\sum_{t=p+1}^{T} \left [\frac{\exp (\bm{Y_{-t}'\beta})}{\{1+\exp (\bm{Y_{-t}'\beta})\}^2}\bm{Y_{-t}}\bm{Y_{-t}}' \right ]\label{heg}.
\end{equation} \qed

\begin{remark}
Note that the results given in (2) and (3) depend on the true values of $\bm{\beta}$. In practice, one needs to plug the maximum likelihood estimates into those expressions to obtain the exact values of Ex-FI.
\end{remark}


\subsection{Logistic autoregressive model of order $p$ with exogenous covariates (LARX(p))}
Here we consider the case of additional exogenous covariate adjustment in the LAR(p) time series model. 
\begin{corollary} 
	\label{cor1}
	Consider a binary-valued correlated time series $Y_t, \ t=1, \ldots, T$, where the conditional distribution of $Y_t$ depends on its previous values and exogenous covariates $X_t$ that relates to current time $t$ through the conditional probability
\[
P_t  = \prob (Y_t=1\mid y_{t-1}, y_{t-2}, \cdots, y_1)=\frac{\exp (\bm{y_{-t}'\beta}+\bm{x_{t}'\alpha})}{1+\exp (\bm{y_{-t}'\beta}+\bm{x_{t}'\alpha})}.
\]
where $\bm{x_{t}}=(x_{t1},\cdots, x_{tl})'$, $\bm{\alpha}=(\alpha_1, \cdots, \alpha_l)'$ and all the other parameters follow the notation of the previous section.
The Ex-FI matrix takes the form
\begin{equation}
I(\bm{\alpha, \beta}\mid y_{p},\cdots, y_1)= \sum_{t=p+1}^{T}\sum\limits_{\bm{(y_{-t})}}\Big{[}\frac{\exp (\bm{X_{t}'\alpha+Y_{-t}'\beta})}{\{1+\exp (\bm{X_{t}'\alpha+Y_{-t}'\beta})\}^2}\left(\begin{array}{cc}
\bm{X_{t}X_{t}'} &\bm{X_{t}Y_{-t}'} \\
\bm{Y_{-t}X_{t}'} & \bm{Y_{-t}Y_{-t}'}
\end{array}\right)\Big{]} Q_t (y_{t-1}, \cdots, y_{t-p}),
\end{equation}
where $Q_t (y_{t-1}, \cdots, y_{t-p})$ is defined in Equation (\ref{marginal}). 
\end{corollary}

\textbf{Proof} The results follow directly from Theorem \ref{thm1} and the facts that the conditional log-likelihood function of $\bm{\alpha}, \bm{\beta}$ and the vector of conditional score functions are respectively
\begin{eqnarray}
\label{likelax}
\ell(\bm{\alpha, \beta}\mid \bm{Y}) & = & \sum_{t=p+1}^{T}[Y_t(\bm{X_{t}'\alpha+Y_{-t}'\beta})-\log\{1+\exp(\bm{X_{t}'\alpha+Y_{-t}'\beta})\}], \\
\bm{U(\alpha, \beta, Y)} & = & \sum_{t=p+1}^{T}\left(\begin{array}{c}
\bm{X_{t}}\{Y_t-\frac{\exp (\bm{X_{t}'\alpha+Y_{-t}'\beta})}{1+\exp (\bm{X_{t}'\alpha+Y_{-t}'\beta})}\}\\
\bm{Y_{-t}}\{Y_t-\frac{\exp (\bm{X_{t}'\alpha+Y_{-t}'\beta})}{1+\exp (\bm{X_{t}'\alpha+Y_{-t}'\beta})}\}
\end{array}\right). \nonumber
\end{eqnarray}
And the Hessian matrix is
\begin{equation}  
\bm{H(\alpha, \beta \mid Y)}= -\sum_{t=p+1}^{T} \left [\frac{\exp (\bm{X_{t}'\alpha+Y_{-t}'\beta})}{\{1+\exp (\bm{X_{t}'\alpha+Y_{-t}'\beta})\}^2}\left(\begin{array}{cc}
\bm{X_{t}X_{t}'} &\bm{X_{t}Y_{-t}'} \\
\bm{Y_{-t}X_{t}'} & \bm{Y_{-t}Y_{-t}'}
\end{array}\right) \right ].\label{heg2}
\end{equation} \qed

In practice, examples of exogenous $X_t$ have been discussed in \citet{Davis:2000}. A particular example is $X_t=t/n$ if one believes that there is a linear temporal trend in the link function (e.g., log mean for the Poisson response in \citet{Davis:2000} and the log odds for LARX(p)).

\subsection{Computation through functional iteration}
Since $(y_{t-1}, \cdots, y_{t-p})$ take values from $\{0,1\}^p$, computation of the Ex-FI matrix through direct calculation can be expensive. In this section, we propose two alternative approaches to achieve Ex-FI matrix. 
\begin{theorem}
	\label{thm2}
	The Ex-FI matrices of Theorem \ref{thm1} and Corollary \ref{cor1} can be achieved by 
the recursive relationship that 

\hskip -1.5cm
	\vbox{\begin{align*}
		& Q_t(y_{t-1}, \cdots, y_{t-p})\\
		&=\left\{
		\begin{array}{ll}
		\sum\limits_{w\in\{0,1\}}\{\prob (Y_{t-1}=1\mid y_{t-2}, \cdots, y_2, w)^{y_{t-1}}\prob (Y_{t-1}=0\mid y_{t-2}, \cdots, y_2, w)^{(1-{y_{t-1}})}& \quad \text{if} \quad p+2\leq t \leq T\\
		\times Q_{t-1}(y_{t-2}, \cdots, y_{t-p}, w)\},&\\
		\mathbbm{1}(y_{p}, \cdots, y_{1}),  & \quad  \text{if} \quad t=p+1. 
		\end{array}
		\right.
		\end{align*}}
where $\mathbbm{1}(y_{p}, \cdots, y_{1})$ is the indicator function that takes value $1$ when the realization is $(y_{p}, \cdots, y_{1})$ and $0$ otherwise.

Alternatively, the Ex-FI matrix can be also obtained through 
	$
	I(\bm{\beta}\mid y_{p},\cdots, y_1)= \sum \limits_{t_0=p+1}^{T}f_{t_0-p}(\bm{\tilde{y}_{-(p+1)}}),
	$ 
	where 
	\begin{align*}
	& f_k(\bm{\tilde{y}_{-(t-k+1)}})\\
	&=\left\{
	\begin{array}{lr}
	f_{k-1}(\bm{\tilde{y}^0_{-(t-k+1)}})+(f_{k-1}(\bm{\tilde{y}^1_{-(t-k+1)}})-f_{k-1}(\bm{\tilde{y}^0_{-(t-k+1)}}))\frac{\exp (\bm{y_{-(t-k+1)}'\beta})}{1+\exp (\bm{y_{-(t-k+1)}'\beta})}, & \text{if} \quad 2\leq k \leq t-p\\
	\frac{\exp (\bm{y_{-t}'\beta})}{\{1+\exp (\bm{y_{-t}'\beta})\}^2}\bm{y_{-t}}\bm{y_{-t}}',  & \text{if} \quad k=1,
	\end{array}
	\right.
	\end{align*}
$\bm{\tilde{y}_{-t}}=(y_{t-1}, \cdots, y_{t-p})', \bm{\tilde{y}^0_{-t}}=(0, y_{t-2}, \cdots, y_{t-p})'$ and $\bm{\tilde{y}^1_{-t}}=(1, y_{t-2}, \cdots, y_{t-p})'.$
\end{theorem}

\textbf{Proof} The derivation of $Q_t(y_{t-1}, \cdots, y_{t-p})$ directly follows from the definition. The results of $f_k(\bm{\tilde{y}_{-(t-k+1)}})$ derive from the fact that 
for any particular $t_0$, by function iteration, 
\begin{align*}
\mathbb{E} \left[ \frac{\exp (\bm{Y_{-t_0}'\beta})}{\{1+\exp (\bm{Y_{-t_0}'\beta})\}^2}\bm{Y_{-t_0}}\bm{Y_{-t_0}}' \right]=f_{t_0-p}(\bm{\tilde{y}_{-(p+1)}}), \quad  t_0 \geq p+1.\\
\end{align*} \qed

\begin{remark}
	From Theorem \ref{thm2}, the marginal probability mass function $Q_t (y_{t-1}, \cdots, y_{t-p})$ can be obtained iteratively at low computational cost. The Ex-FI matrix can be achieved accordingly. On the other hand, the Ex-FI matrix can also be obtained via iterated expectations from the second part of Theorem \ref{thm2}.
\end{remark}

\begin{remark}
	The results from Theorem~\ref{thm2} can be extended to other link functions. In the case of complementary log-log link, given by the analytic form of the latent function, score function and Hessian matrix can be easily obtained. Equations~(4) and (7) can be directly adapted. And by following Theorem 2, the Ex-FI can be obtained. In the case of probit link, due to the inexplicit form of the link function, a numerical approximation can be used to obtain Ex-FI in practice.
\end{remark}


\subsection{Special case: logistic autoregressive model of order $p=1$ (LAR(1))}
In general, there is no explicit form of Ex-FI. In this section, we consider the only special case that enjoys explicit analytic form.
\begin{theorem}
Consider a binary-valued time series data $Y_t, \ t=1, \ldots, T$, where the conditional distribution of $Y_t$ depends on its
own immediate past value via the conditional probability
\[
P_t  = \prob (Y_t=1\mid y_{t-1}, y_{t-2}, \cdots, y_1)=\frac{\exp (\beta_0+\beta_1 y_{t-1})}{1+\exp (\beta_0+\beta_1 y_{t-1})}.
\]
Then if we denote the Ex-FI matrix to be $I(\bm{\beta}\mid Y_1)$, its elements $I_{jk}, j=1,2; k=1,2$ are derived,
respectively, as
\begin{eqnarray*}
	I_{11} & = & \mathbb{E}[\sum\limits_{t=3}^{T} v(Y_{t-1}) \mid Y_1]+v(Y_1)\\
	& = & \sum\limits_{t=3}^{T} [\{v(1)-v(0)\}\{{p(1)}-{p(0)}\}^{t-3}\{p(Y_1)-\frac{p(0)}{1-p(1)+p(0)}\} + \\
	{} & {} &  v(0)+p(0)\frac{v(1)-v(0)}{1-p(1)+p(0)}]+v(Y_1)\\
	{} &= & \{v(1)-v(0)\}\{p(Y_1)-\frac{p(0)}{1-p(1)+p(0)}\}\frac{1-\{p(1)-p(0)\}^{T-2}}{1-p(1)+p(0)}+\\
	{} & {} & \frac{(T-2)\{p(0)v(1)+v(0)-v(0)p(1)\}}{1-p(1)+p(0)}+v(Y_1)
\end{eqnarray*}
\begin{eqnarray*}
	I_{12} & = & \mathbb{E}[\sum\limits_{t=3}^{T}  v(Y_{t-1}) Y_{t-1} \mid Y_1]+v(Y_1)Y_1\\
	{} & = & \sum\limits_{t=3}^{T} [v(1)\{{p(1)}-{p(0)}\}^{t-3}\{p(Y_1)-\frac{p(0)}{1-p(1)+p(0)}\}+\frac{v(1)p(0)}{1-p(1)+p(0)}]+ v(Y_1)Y_1 \\
	{} & = & v(1)\{p(Y_1)-\frac{p(0)}{1-p(1)+p(0)}\}\frac{1-\{p(1)-p(0)\}^{T-2}}{1-p(1)+p(0)}+\frac{(T-2)p(0)v(1)}{1-p(1)+p(0)}+v(Y_1)Y_1\\
	I_{22} & = & \mathbb{E}[\sum\limits_{t=3}^{T} v(Y_{t-1}) Y_{t-1}^2  \mid Y_1]+v(Y_1)Y_1^2\\
	&=&I_{12},
\end{eqnarray*}
where  $p(y)=\frac{\exp(\beta_0+\beta_1y)}{ 1+\exp(\beta_0+\beta_1y)}$ and $v(y) = p(y)[1-p(y)]$.
\end{theorem}

\textbf{Proof}
It is straightforward that the corresponding Hessian matrix is derived to be
\begin{eqnarray} \label{he}
\frac{\partial ^2}{\partial \bm{\beta} \partial \bm{\beta^T}} \ell(\bm{\beta}\mid \bm{Y})= - \left(
\begin{array}{cc}
\sum\limits_{t=2}^{T} v(Y_{t-1})  & \sum\limits_{t=2}^{T} v(Y_{t-1}) Y_{t-1} \\
\sum\limits_{t=2}^{T} v(Y_{t-1}) Y_{t-1} & \sum\limits_{t=2}^{T} v(Y_{t-1}) Y_{t-1}^2 \\
\end{array}
\right).
\end{eqnarray}
Due to the Markovian assumption, the conditional expectation can be obtained
through iterated expectations.
For any particular $t>2$, we have
\begin{eqnarray}
\mathbb{E} \left [ v(Y_{t-1}) \mid Y_1 \right] & = &  A_1 + A_2 \label{e1} \\
\mathbb{E} \left [ v(Y_{t-1}) Y_{t-1} \mid Y_1 \right] & = & A_3 + A_4  \ \ {\mbox{where}} \label{e2}
\end{eqnarray}
\begin{align*}
A_1 & =  \{v(1)-v(0)\} \{p(1)-p(0)\}^{t-3}\{p(Y_1)- p(0)/\{1-p(1)+p(0)\} \}\\
A_2 & =  v(0)+ p(0)\{v(1)-v(0)\}/\{1-p(1)+p(0)\} \\
A_3 & =  v(1)\{p(1)-p(0)\}^{t-3} \left[ p(Y_1)- p(0)/\{1-p(1)+p(0)\} \right] \\
A_4 & =  v(1)p(0)/\{1-p(1)+p(0)\},
\end{align*}
which completes the proof by some algebra calculation.  \qed

\section{Simulations}
\label{section:sim}
\subsection{Evaluating small sample performance}
In this section, we compare the behavior of the newly derived Ex-FI and Em-FI in the context of inference for regression parameters under various models. The Em-FI (AFI) is calculated from Equation (4) for LAR and Equation (7) for LARX respectively. Time series lengths are chosen as $T=20, 50$ and $200$ respectively, and $10,000$ simulations were generated under different scenarios. In Scenario 1, the signals were generated by LAR(1) with parameters $(\beta_0, \beta_1)=(0.1, 0.5)$ (``low ratio"), and $(\beta_0, \beta_1)=(0.1, 1)$ (``high ratio"). In this case, $\beta_1$ denotes the log odds ratio and  $\beta_0$ denotes the log odds when the previous realization is $0$. $\beta_1/\beta_0$ is a monotonic function of the log odds ratio of $Y_t=1$. Particularly, large value (greater than 1) of $\beta_1/\beta_0$ implies the log odds of $Y_t=1$ when $Y_{t-1}=1$ is much higher compared to the log odds when $Y_{t-1}=0$. In Scenario 2, the time series were simulated through LAR(2) with parameters  $(\beta_0, \beta_1, \beta_2)=(0.1, 0.3, 0.5)$ (``low ratio"), and $(\beta_0, \beta_1, \beta_2)=(0.1, 1, 1.5)$ (``high ratio"). In Scenario 3, we considered generating signals by LARX(1) with parameters $(\alpha_1, \beta_0, \beta_1)=(0.5, 0.1, 0.5)$ (``low ratio"), and $(\alpha_1, \beta_0, \beta_1)=(0.5, 0.1, 1)$ (``high ratio"). The exogenous covariate was obtained from standard normal distribution. For each scenario, we calculate the empirical type I error-rate for testing $H_0:\beta_1 (\alpha_1, \beta_2)=0$ at level .05, the average standard error obtained from maximum likelihood estimates and true values as well as their corresponding Monte Carlo standard error, the observed error deviation of the estimates across simulations. Critical values are determined from normal distribution when calculating type I error-rates.

Table \ref{table1} provides a summary of the conducted simulation study for various time series lengths.  With respect to type I error, it can be seen that use of Ex-FI and Em-FI both result in conservative inference (lower than nominal type I error) for smaller values of $T$ and for high ratios.  For the low ratio scenario, nominal type I error rates are achieved as time series lengths of $T=50$. For time series lengths of $T=200$ both variance estimators yield the desired type I error rates. 
The benefit of using Ex-FI over Em-FI is observed when comparing the average standard error to the observed standard deviation of estimates of $\beta_1$ across simulations.  Specifically, Em-FI tends to behave erratically for small sample sizes, yielding extremely large estimated standard error for some simulated datasets.  This can be seen most notably in the high ratio scenario by observing that the average standard error computed using Em-FI is 362.3 compared to the actual observed standard deviation of the estimator across 10,000 simulations being only 3.015.  In contrast, the average standard error computed using Ex-FI is only 7.868. 

Table~\ref{table2} summarizes the simulation results of Scenario 2. Similar to the case of LAR(1), both of the Ex-FI and Em-FI result in lower type I error rates when the sample size is small and reach nominal values when $T = 200.$ We can still observe the huge advantages of the average standard error obtained from Ex-FI in contrast to the ones from Em-FI. From the Monte Caro errors presented by the second values in the parenthesis, such advantages are statistical significant in most of the cases especially when $T$ is relatively small. Moreover, by comparing between the first value in the parenthesis to the observed standard error, one can easily find that the proposed Ex-FI (by inserting the true values of parameters) is close to the Monte Carlo standard error of maximum likelihood estimates. Table~\ref{table2b} presents the results of Scenario 3. Similar findings can be easily found as well.

As is shown in Table~\ref{table1}, in the cases of $T=20, 50$, although the average Ex-FI obtained from maximum likelihood estimates is not very close to the observed standard error of the estimates, the ones obtained from the true values of $\beta_1$ are much closer. In the worst case of $T=20$ (low ratio), the difference is $0.03$ for Ex-FI in comparison with $10$ for Em-FI. Such advantage is much more obvious in the scenario of high ratio. As $T$ increases, the Ex-FI obtained from true values is converging to the observed standard error. From Tables~\ref{table2} and \ref{table2b}, we can clearly find the same pattern for all the parameters of endogenous and exogenous covariates. When the sample size $T$ is relatively small, the Ex-FI obtained from the true value of parameters is close to the observed standard error. The discrepancy is approaching to 0 as $T$ goes towards 200.

\begin{table}[H]
	\caption{Summary of simulation results for the LAR(1) model.  Time series lengths are chosen as $T=20, 50$ and $200$. $10,000$ simulations were generated under each of two parameters settings: $(\beta_0, \beta_1)=(0.1, 0.5)$ (``low ratio") and $(\beta_0, \beta_1)=(0.1, 1)$ (``high ratio").  For each scenario, we present the empirical type I error-rate for testing $H_0: \beta_1=0$, the (average) standard error of the point estimate (and true value) for $\beta_1$, the Monte Carlo standard error, the observed standard error of the regression parameter estimate of $\beta_1$ across simulations.}
	\hskip -5em
	\begin{scriptsize}
		\begin{tabular}{llccccccc}
			\hline
			&& \multicolumn{3}{c}{Low Ratio ($(\beta_0, \beta_1)=(0.1, 0.5)$)} && \multicolumn{3}{c}{High Ratio ($(\beta_0, \beta_1)=(0.1, 1)$)}\\
			\cline{3-5} \cline{7-9}
			\multicolumn{2}{c}{Length/Method} & Type I   & Standard & Observed Standard &&  Type I & Standard & Observed Standard \\ 
			&& Error  & Error* & Error** &&  Error  & Error* & Error**\\ 
			\hline
			$T=20$  \\
			& Ex-FI	&	0.031	&  $\bm{3.737} (2.320, 0.324)$   & \multirow{2}{*}{2.290}
			&&	0.008        & $\bm{7.868} (3.075, 0.543)$ 	&   \multirow{2}{*}{3.015} \\
			& Em-FI   	&      0.030        &  $\bm{32.870} (12.290, 2.438)$   & 	&&    0.011 	& $\bm{362.300} (30.343, 14.355)$ 	& \\
			$T=50$\\  
			& Ex-FI    	&      0.048   	&   $\bm{0.617} (0.630, 0.084)$ 	&  \multirow{2}{*}{0.632}   	&&    0.039       	& $\bm{1.065} (1.070, 0.064)$ 	&  \multirow{2}{*}{1.074} \\
			& Em-FI      &    	0.044        &   $\bm{0.956} (0.748, 0.055)$ 	&   	&&    0.039       	& $\bm{1.222} (1.148, 0.056)$ 	& \\
			$T=200$\\
			& Ex-FI       &      0.052   	&   $\bm{0.299} (0.299, 0.006)$ 	&  \multirow{2}{*}{0.299}   	&&  	 0.051   	& $\bm{0.332} (0.326, 0.008)$	&  \multirow{2}{*}{0.325} \\
			& Em-FI     	&      0.052     	&   $\bm{0.297} (0.298, 0.004)$ 	&     &&     0.053     	& $\bm{0.324} (0.325, 0.005)$	&  \\ 
			\hline
			\multicolumn{9}{l}{*Standard error represents the average standard error of the point estimator for $\beta_1$ and, in parentheses, the same value obtained from true $\beta_1$}\\
			\multicolumn{9}{l}{~~ followed by the Monte Carlo standard error. Note that since Ex-FI does not rely on the true realizations $Y_t$, it remains the same within }\\
			\multicolumn{9}{l}{ ~~ the same scenario while for Em-FI, the reported value is average over repetitions. }\\
			\multicolumn{9}{l}{**Observed standard error represents the Monte Carlo standard error of the maximum likelihood estimates of $\beta_1$ across simulations.}
		\end{tabular}
	\end{scriptsize}
	\label{table1}
\end{table} 

\begin{table}[H]
	\caption{Summary of simulation results for the LAR(2) model.  Time series lengths are chosen as $T=20, 50$ and $200$. $10,000$ simulations were generated under each of two parameters settings: $(\beta_0, \beta_1, \beta_2)=(0.1, 0.3, 0.5)$ (``low ratio") and $(\beta_0, \beta_1, \beta_2)=(0.1, 1, 1.5)$ (``high ratio").  For each scenario, we present the empirical type I error-rate for testing $H_0: \beta_1 (\beta_2) =0$, the (average) standard error of the point estimator (and true value), the Monte Carlo standard error, the observed standard error of the regression parameter estimates across simulations.}
	\hskip -4em
	\begin{scriptsize}
		\begin{tabular}{llccccccc}
			\hline
			&& \multicolumn{3}{c}{Low Ratio ($(\beta_0, \beta_1, \beta_2)=(0.1, 0.3, 0.5)$)} && \multicolumn{3}{c}{High Ratio ($(\beta_0, \beta_1, \beta_2)=(0.1, 1, 1.5)$)}\\
			\cline{3-5} \cline{7-9}
			\multicolumn{2}{c}{Length/Method} & Type I   & Standard & Observed Standard &&  Type I & Standard & Observed Standard \\ 
			&& Error  & Error* & Error** &&  Error  & Error* & Error**\\ 
			\hline
			$T=20$  \\
			\multirow{2}{*}{$\beta_1$}	& Ex-FI	&	0.030	&  $\bm{8.960} (4.351, 2.354)$   & \multirow{2}{*}{4.312}
			&&	0.027        & $\bm{10.149} (8.102, 2.895)$ 	&   \multirow{2}{*}{7.944} \\
			& Em-FI   	&      0.031        &  $\bm{294.409} (56.294, 21.343)$   & 	&&    0.028 	& $\bm{254.576} (42.135, 20.540)$ 	& \\
			\multirow{2}{*}{$\beta_2$}	& Ex-FI	&	0.042	&  $\bm{8.001} (3.942, 2.103)$   & \multirow{2}{*}{3.901}
			&&	0.028        & $\bm{7.868} (7.041, 3.012)$ 	&   \multirow{2}{*}{6.931} \\
			& Em-FI   	&      0.041        &  $\bm{363.611} (64.239, 24.031)$   & 	&&    0.031 	& $\bm{190.567} (36.356, 27.012)$ 	& \\

			$T=50$  \\
			\multirow{2}{*}{$\beta_1$}	& Ex-FI	&	0.048	&  $\bm{1.031} (1.250, 0.073)$   & \multirow{2}{*}{1.247}
			&&	0.030        & $\bm{3.339} (3.286, 0.054)$ 	&   \multirow{2}{*}{3.284} \\
			& Em-FI   	&      0.047        &  $\bm{1.544} (1.532, 0.085)$   & 	&&    0.031 	& $\bm{6.433} (1.845, 0.125)$ 	& \\
			\multirow{2}{*}{$\beta_2$}	& Ex-FI	&	0.042	&  $\bm{1.011} (0.942, 0.051)$   & \multirow{2}{*}{0.949}
			&&	0.043        & $\bm{5.025} (3.995, 0.083)$ 	&   \multirow{2}{*}{3.993} \\
			& Em-FI   	&      0.041        &  $\bm{1.836} (1.825, 0.044)$   & 	&&    0.042 	& $\bm{5.806} (4.024, 0.121)$ 	& \\

			$T=200$  \\
			\multirow{2}{*}{$\beta_1$}	& Ex-FI	&	0.052	&  $\bm{0.734} (0.614, 0.008)$   & \multirow{2}{*}{0.612}
			&&	0.047        & $\bm{2.562} (2.521, 1.042)$ 	&   \multirow{2}{*}{2.522} \\
			& Em-FI   	&      0.051        &  $\bm{0.849} (0.753, 0.004)$   & 	&&    0.045 	& $\bm{4.834} (3.454, 1.021)$ 	& \\
			\multirow{2}{*}{$\beta_2$}	& Ex-FI	&	0.048	&  $\bm{0.801} (0.702, 0.007)$   & \multirow{2}{*}{0.701}
			&&	0.048        & $\bm{1.762} (1.504, 0.542)$ 	&   \multirow{2}{*}{1.503} \\
			& Em-FI   	&      0.050        &  $\bm{0.913} (0.645,0.004)$   & 	&&    0.050 	& $\bm{1.864} (1.735, 0.842)$ 	& \\
			\hline
			\multicolumn{9}{l}{*Standard error represents the average standard error of the point estimates and, in parentheses, the same measure obtained from true values}\\
			\multicolumn{9}{l}{~~ followed by the Monte Carlo standard error. Note that since Ex-FI does not rely on the true realizations $Y_t$, it remains the same within }\\
			\multicolumn{9}{l}{ ~~ the same scenario while for Em-FI, the reported value is average over repetitions. }\\
			\multicolumn{9}{l}{**Observed standard error represents the Monte Carlo standard error of the maximum likelihood estimates across simulations.}
		\end{tabular}
	\end{scriptsize}
	\label{table2}
\end{table} 

\begin{table}[H]
	\caption{Summary of simulation results for the LARX(1) model.  Time series lengths are chosen as $T=20, 50$ and $200$. $10,000$ simulations were generated under each of two parameters settings: $(\alpha_1, \beta_0, \beta_1)=(0.5, 0.1, 0.5)$ (``low ratio") and $(\alpha_1, \beta_0, \beta_1)=(0.5, 0.1, 1)$ (``high ratio").  For each scenario, we present the empirical type I error-rate for testing $H_0: \beta_1 (\alpha_1)=0$, the (average) standard error of the point estimator (and true value), the Monte Carlo standard error, the observed standard error of the regression parameter estimate  across simulations.}
	\hskip -4em
	\begin{scriptsize}
		\begin{tabular}{llccccccc}
			\hline
			&& \multicolumn{3}{c}{Low Ratio ($(\alpha_1, \beta_0, \beta_1)=(0.5, 0.1, 0.5)$)} && \multicolumn{3}{c}{High Ratio ($(\alpha_1, \beta_0, \beta_1)=(1, 0.1, 1)$)}\\
			\cline{3-5} \cline{7-9}
			\multicolumn{2}{c}{Length/Method} & Type I   & Standard & Observed Standard &&  Type I & Standard & Observed Standard \\ 
			&& Error  & Error* & Error** &&  Error  & Error* & Error**\\ 
			\hline
			$T=20$  \\
			\multirow{2}{*}{$\alpha_1$}	& Ex-FI	&	0.027	&  $\bm{10.801} (6.382, 2.753)$   & \multirow{2}{*}{6.363}
			&&	0.031        & $\bm{18.535} (17.302, 4.435)$ 	&   \multirow{2}{*}{17.522} \\
			& Em-FI   	&      0.029        &  $\bm{241.515} (24.352, 8.954)$   & 	&&    0.032 	& $\bm{352.153} (31.233, 14.983)$ 	& \\
			\multirow{2}{*}{$\beta_1$}	& Ex-FI	&	0.032	&  $\bm{13.242} (5.942, 2.021)$   & \multirow{2}{*}{6.018}
			&&	0.038        & $\bm{17.322} (12.011, 4.321)$ 	&   \multirow{2}{*}{11.460} \\
			& Em-FI   	&      0.035        &  $\bm{134.542} (34.240, 10.324)$   & 	&&    0.036	& $\bm{179.222} (14.324, 9.921)$ 	& \\

			$T=50$  \\
			\multirow{2}{*}{$\alpha_1$}	& Ex-FI	&	0.048	&  $\bm{0.334} (0.332, 0.062)$   & \multirow{2}{*}{0.332}
			&&	0.033        & $\bm{0.566} (0.529, 0.041)$ 	&   \multirow{2}{*}{0.531} \\
			& Em-FI   	&      0.047        &  $\bm{0.852} (0.344, 0.053)$   & 	&&    0.033 	& $\bm{0.963} (0.552, 0.063)$ 	& \\
			\multirow{2}{*}{$\beta_1$}	& Ex-FI	&	0.042	&  $\bm{0.783} (0.694, 0.073)$   & \multirow{2}{*}{0.691}
			&&	0.038        & $\bm{0.785} (0.763, 0.050)$ 	&   \multirow{2}{*}{0.769} \\
			& Em-FI   	&      0.041        &  $\bm{1.333} (0.723, 0.042)$   & 	&&    0.039 	& $\bm{1.420} (0.774, 0.041)$ 	& \\

			$T=200$  \\
			\multirow{2}{*}{$\alpha_1$}	& Ex-FI	&	0.051	&  $\bm{0.194} (0.184, 0.006)$   & \multirow{2}{*}{0.185}
			&&	0.048        & $\bm{0.198} (0.198, 0.003)$ 	&   \multirow{2}{*}{0.198} \\
			& Em-FI   	&      0.050        &  $\bm{0.199} (0.193, 0.007)$   & 	&&    0.046 	& $\bm{0.196} (0.196, 0.004)$ 	& \\
			\multirow{2}{*}{$\beta_1$}	& Ex-FI	&	0.049	&  $\bm{0.315} (0.314, 0.007)$   & \multirow{2}{*}{0.314}
			&&	0.050        & $\bm{0.343} (0.343, 0.005)$ 	&   \multirow{2}{*}{0.343} \\
			& Em-FI   	&      0.051        &  $\bm{0.316} (0.315, 0.006)$   & 	&&    0.051 	& $\bm{0.342} (0.340, 0.005)$ 	& \\
			\hline
			\multicolumn{9}{l}{*Standard error represents the average standard error of the point estimates and, in parentheses, the same measure obtained from true values}\\
			\multicolumn{9}{l}{~~ followed by the Monte Carlo standard error. Note that since Ex-FI does not rely on the true realizations $Y_t$, it remains the same within }\\
			\multicolumn{9}{l}{ ~~ the same scenario while for Em-FI, the reported value is average over repetitions. }\\
			\multicolumn{9}{l}{**Observed standard error represents the Monte Carlo standard error of the maximum likelihood estimates across simulations.}
		\end{tabular}
	\end{scriptsize}
	\label{table2b}
\end{table} 
\subsection{Evaluation of confidence interval length}
\label{evaci}
\begin{figure}[H] \centering
	\begin{tabular}{cc}
		\includegraphics[width=.5\textwidth]{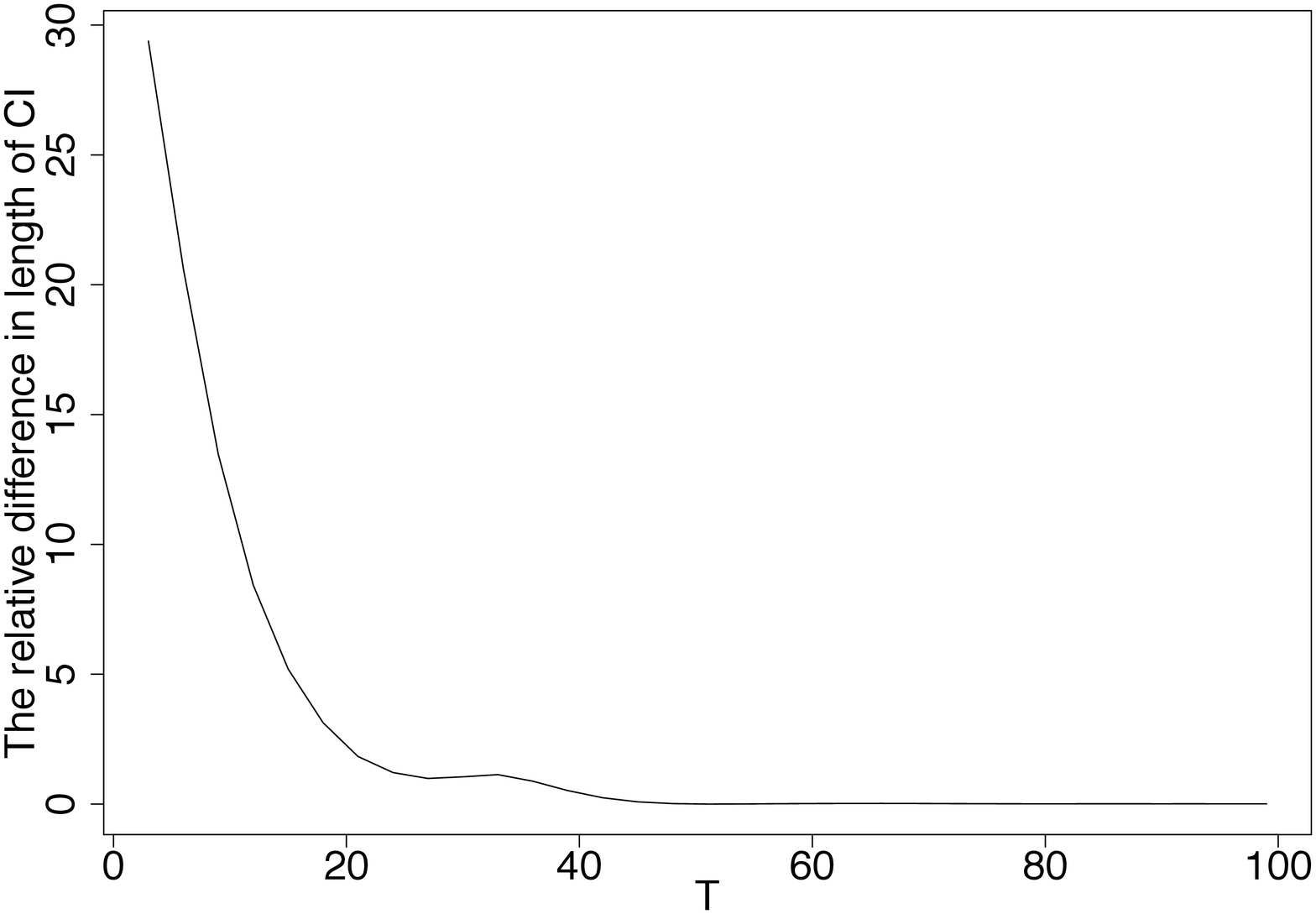} &
		\includegraphics[width=.5\textwidth]{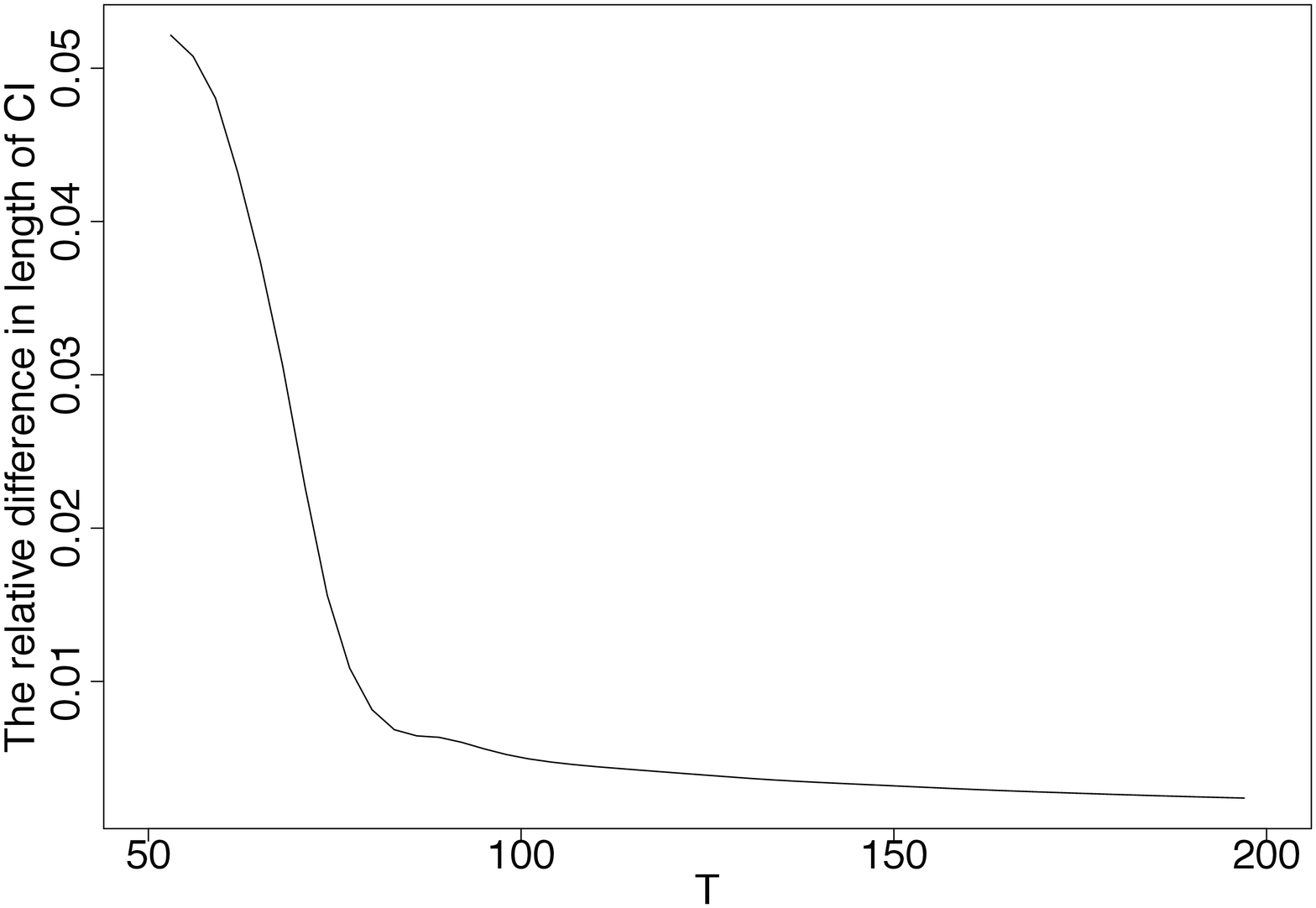} \\
	\end{tabular}
	\caption{The average relative difference in lengths of confidence intervals derived from
		Ex-FI and Em-FI  (length of CI$^{\text{empirical}}-$length of CI$^{\text{exact}}$)/length of CI$^{\text{exact}}$
		computed from $1000$ simulated time series with $\beta_1/{\beta_0}=10$. The lengths of time series, $T$ ranges from (a) 5 to 100 (left) and (b) 50 to 200 (right). The Em-FI matrix used here was identical to the one proposed in \citet{Fokianos:1998}.}
	\label{ci}
\end{figure}

Here we consider the average length of derived 95\% confidence intervals for $\beta_1$. Following the result that asymptotically, $\widehat{\bm{\beta}} \overset{.}{\sim} N(\bm{\beta}, I^{-1}(\bm{\beta}))$ for large values of $T$ \citep{Fokianos:1998}, an approximate $95\%$ confidence interval can be obtained using both Ex-FI and Em-FI. For each scenario of $\bm{\beta}$ described above, 1000 binary time series of lengths $T=5,6, \dots, 250$ were generated. For each time series data, an approximate $95\%$ confidence interval for $\beta_1$ was
computed using both Ex-FI and Em-FI. We compared the two approaches by calculating the relative
difference of the lengths of the two confidence intervals.  
As expected from the average standard error values in Table 1, Fig.~\ref{ci} indicates that the confidence interval derived from Ex-FI behaves more efficiently on average
than the confidence interval computed using Em-FI. It is noted that such substantial difference exists when $T<200$ and tends to be roughly the same as $T$ goes beyond 200. Once again, it implies that one should be careful with the Em-FI when $T<200$.

\begin{figure}[H] \centering
	\begin{tabular}{cc}
		\includegraphics[width=.5\textwidth]{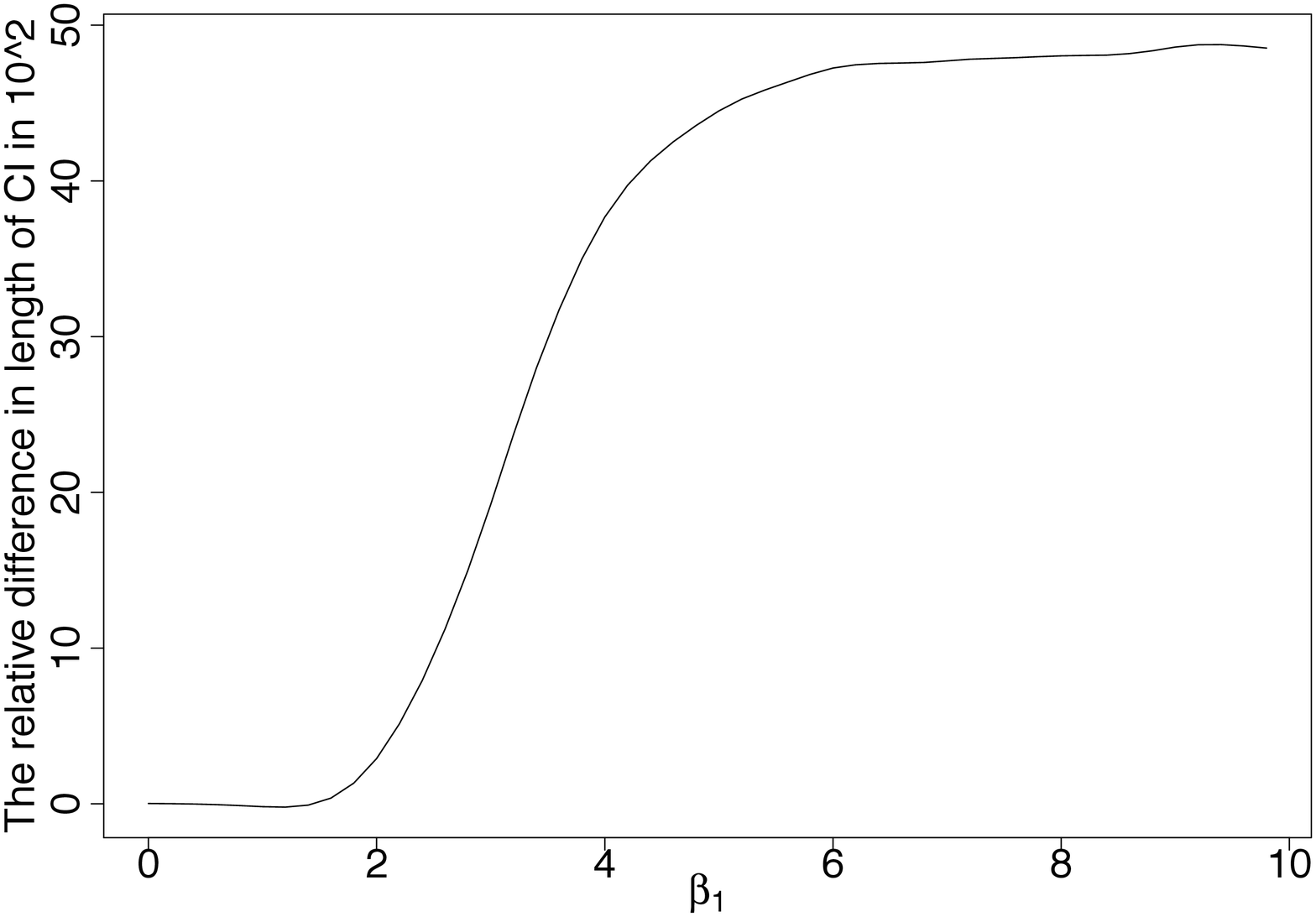} &
		\includegraphics[width=.5\textwidth]{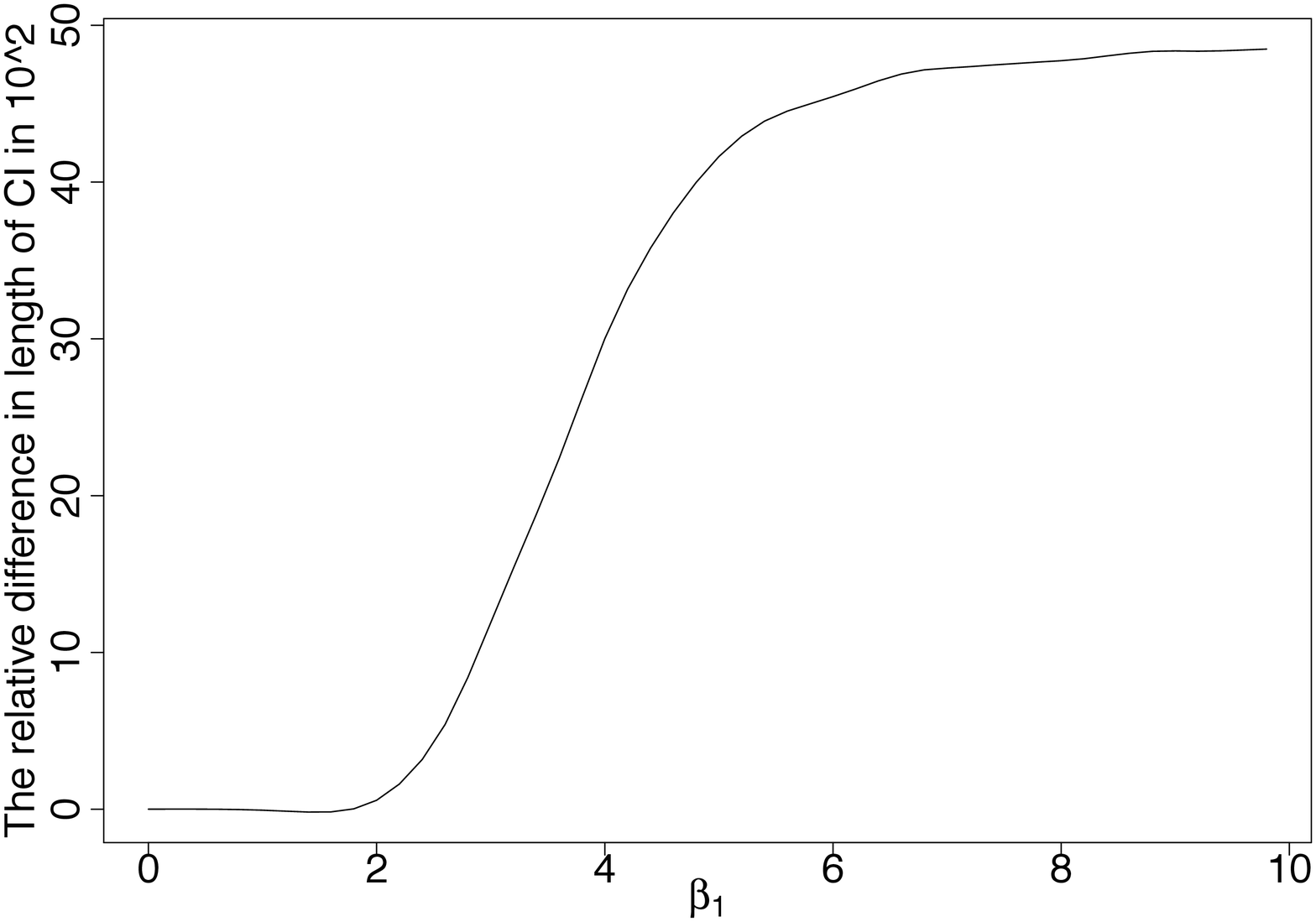} \\
	\end{tabular}
	\caption{The average relative difference in length of confidence intervals, computed
		from $1000$ simulated datasets, derived from
		Ex-FI and Em-FI (length of CI$^{\text{empirical}}-$length of CI$^{\text{exact}}$)/length of CI$^{\text{exact}}$, where the number of observations is taken to be (a) $T=60$ (left) and (b) $T=100$ (right). $\beta_0$ is fixed to be $0.1$. The Em-FI matrix used was developed in \citet{Fokianos:1998}.}
	\label{ci2}
\end{figure}

In Fig.~\ref{ci2}, $T$ was fixed to $60$ and $100$ and $\beta_1$ was allowed to vary while
keeping $\beta_0 = 0.1$. Results clearly establish the advantage of
Ex-FI over Em-FI especially as the true value of $\beta_1$ increases, i.e., the ratio increases.

\subsection{Evaluating the discrepancy between the exact and empirical Fisher information}
\begin{figure}[h] \centering
	\begin{tabular}{cc}
		\includegraphics[width=.5\textwidth]{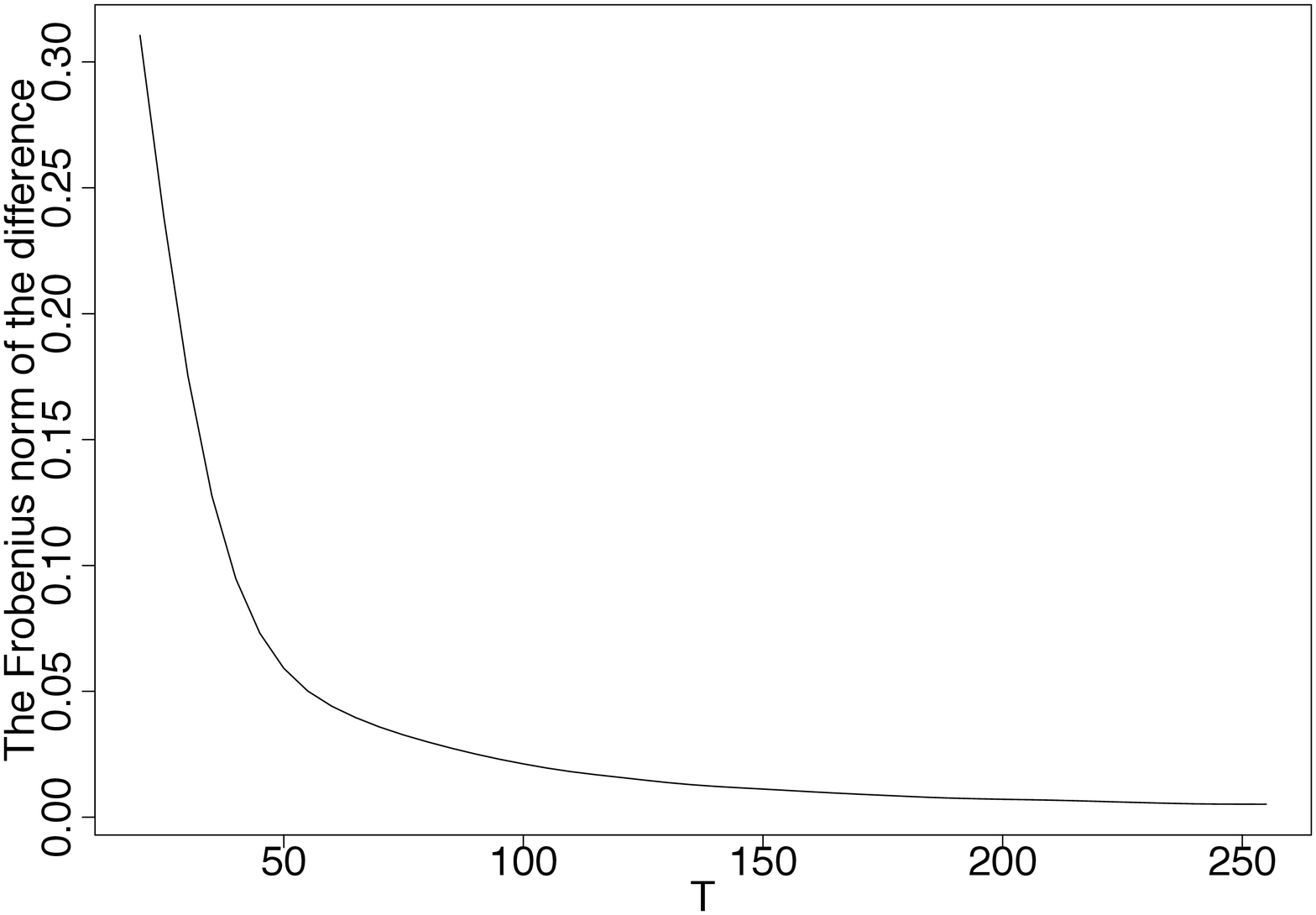} &
		\includegraphics[width=.5\textwidth]{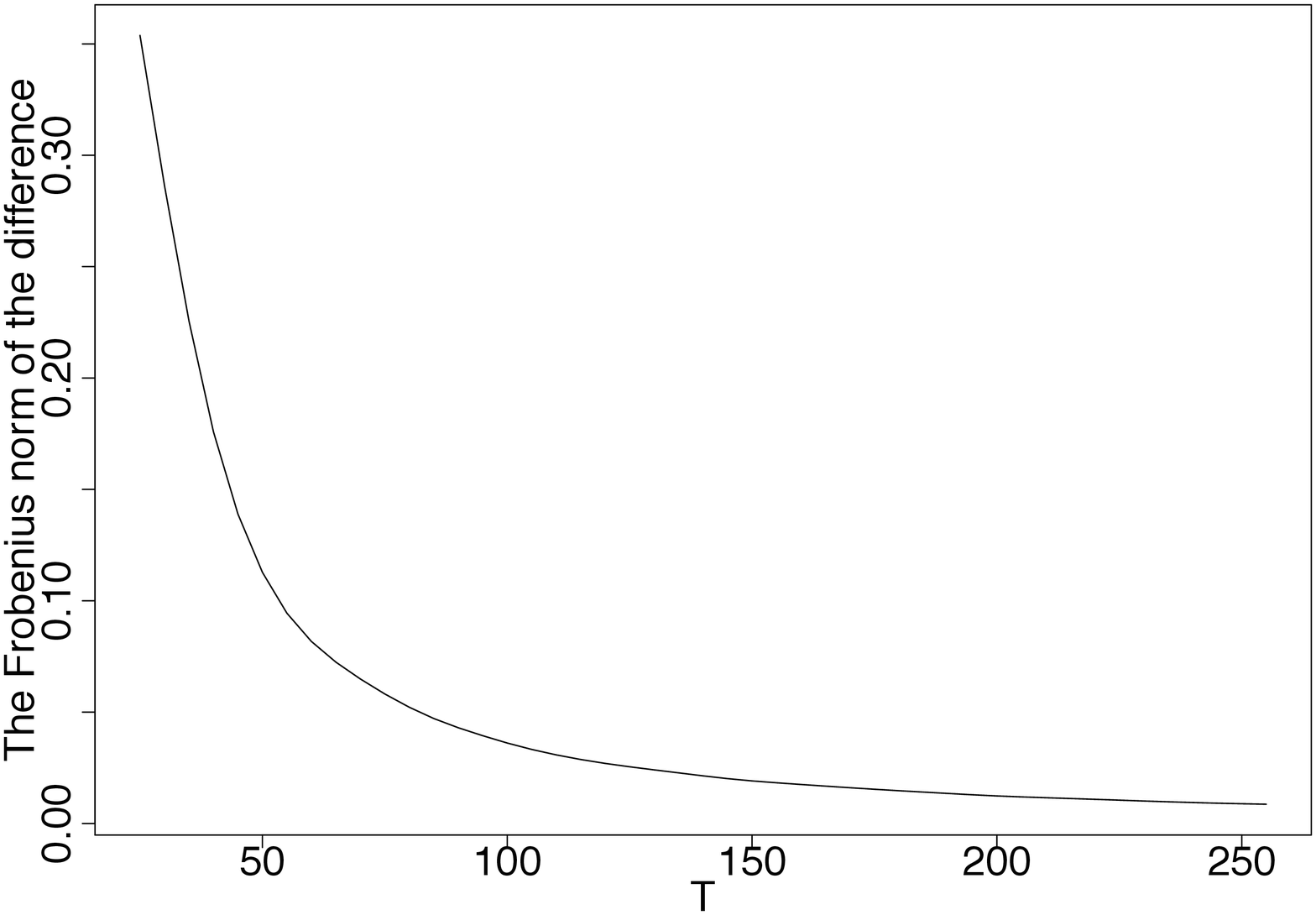} \\
	\end{tabular}
	\caption{The average Frobenius norm of the difference between the inverse of the exact Fisher information (Ex-FI)
		and empirical Fisher information (Em-FI) (as developed in \citet{Fokianos:1998}) under the two parameter set up:  (a) $\beta_1/{\beta_0}=5$ (left) and (b) $\beta_1/{\beta_0}=10$ (right).
		The average Frobenius norm was calculated from 1,000 simulated time series for varying time series lengths under each of the parameter set-up.}
	\label{mse1}
\end{figure}

In this section, we discuss the results of simulations conducted to investigate the discrepancy between Ex-FI and Em-FI under
the following scenarios: (i.) time series lengths $T$ ranging from $10-250$; (ii.) the ratio
$\beta_0/\beta_1$ $\in$ $\{ 5,10 \}$. Based on $1,000$ simulated time series under each scenario, the average Frobenius norm of the difference between the asymptotic covariance matrices (i.e. the inverse of Ex-FI and Em-FI), displayed in Fig.~\ref{mse1}, shows that when $T>200$ any discrepancy between the two covariance matrices effectively vanishes. However, for $T<200$, discrepancies do exist, primarily due to the instability of Em-FI for particular datasets.  The result reiterates that caution needs to be taken when utilizing the Em-FI variance estimator for shorter time series, since this erratic behaviour could lead to significant errors in the estimated variances of regression parameters.

\subsection{Evaluating the convergence}
\begin{figure}[H] \centering
	\begin{tabular}{cc}
		\includegraphics[width=.5\textwidth]{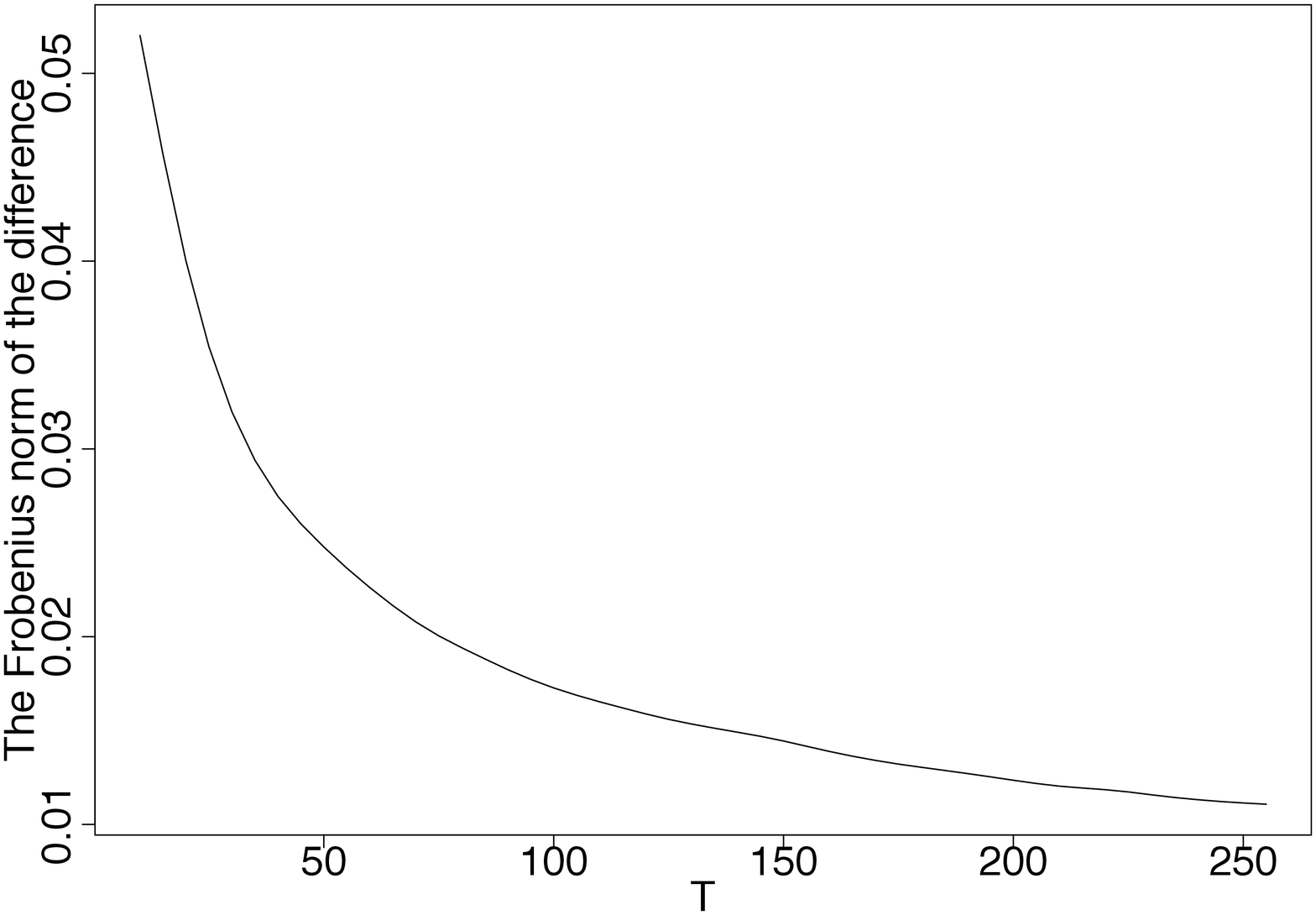} &
		\includegraphics[width=.5\textwidth]{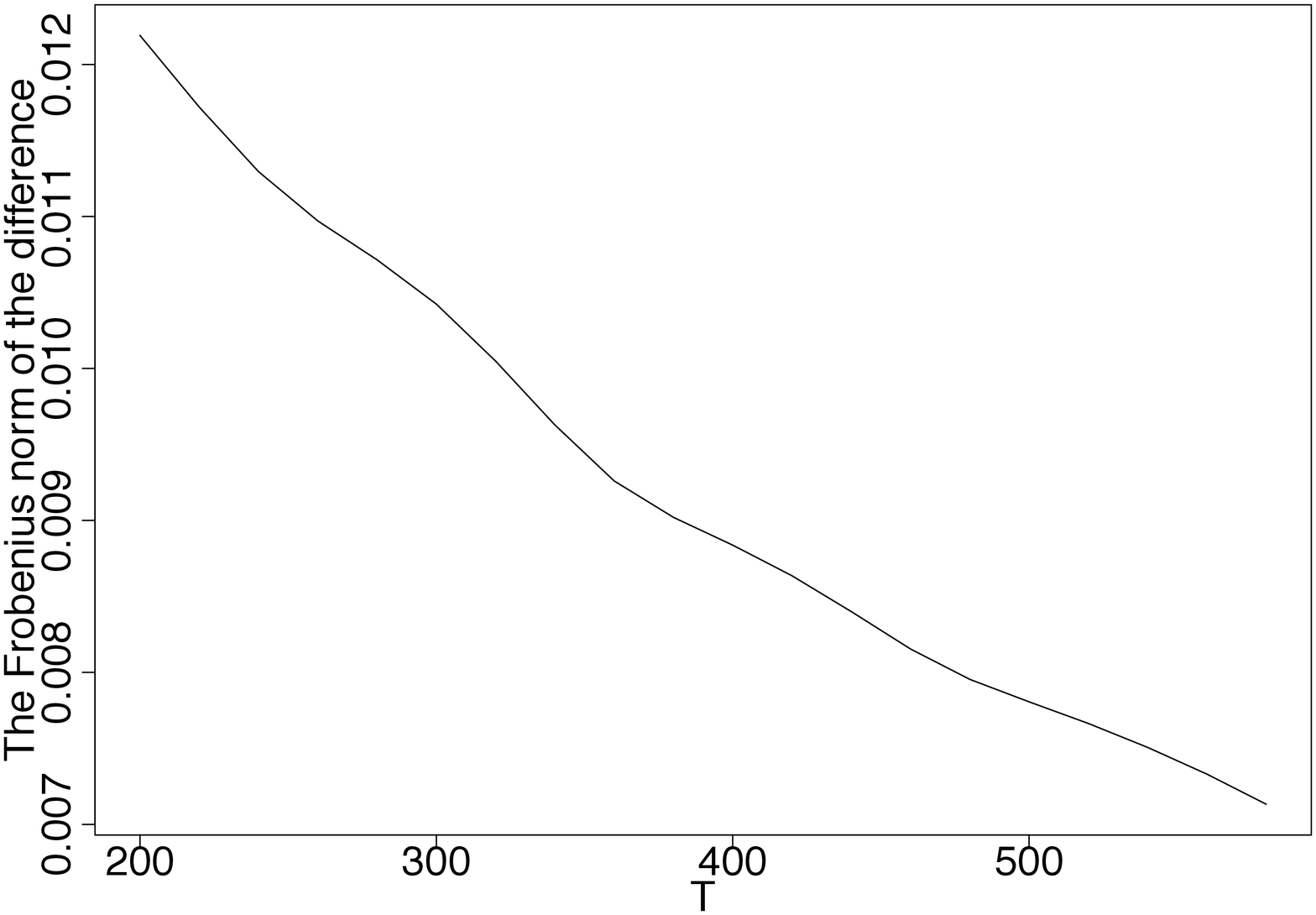} \\
	\end{tabular}
	\caption{The average Frobenius norm of the difference in Ex-FI and AFI matrices (which is proposed in \citet{Fokianos:1998}) computed over $1000$ simulated time series under the set-up $\beta_1/{\beta_0}=5$. The lengths of time series, (a) $T$ ranges from 5 to 250 (left) and (b) 250 to 550 (right).}
	\label{ke1}
\end{figure}

We considered the asymptotic behavior of Ex-FI and compared it to the AFI proposed by \citet{Fokianos:1998}
by computing the average Frobenius norm between the two matrices over $1,000$ simulated time series data.
In Fig.~\ref{ke1}, it is clear that the discrepancy between these two matrices decays dramatically, which empirically indicates that the limiting behavior between the two estimators coincides. It should be emphasized that when $T<200$, the difference is significant while as $T$ grows larger than 200, the discrepancy shrinks to small values around 0. Hence, utilizing the Em-FI when $T<200$ may be problematic.

\section{Analysis of Binary Respiratory Time Series}
\label{section:real}
\subsection{Explanatory analysis}
In this section we consider time-series data on respiratory rate among a cohort of 113 expectant mothers. Briefly, the participants consist of a sub-sample of women from a larger cohort of women attending prenatal care at a university-based clinic in Pittsburgh, PA and participating in a prospective, longitudinal study from early gestation through birth \citep{Entringer:2015}. Participants were asked to wear a heart and respiratory rate monitor for up to four consecutive days.  In addition, each night prior to sleeping the participants were asked to fill out an electronic diary recording how stressful their day was on a scale from 1 to 10 ($X_i$), with 10 corresponding to the highest self-reported stress level. The study was approved by the local Institutional Review Board (IRB). 

\begin{figure}[h] \centering
	\includegraphics[width=.5\textwidth]{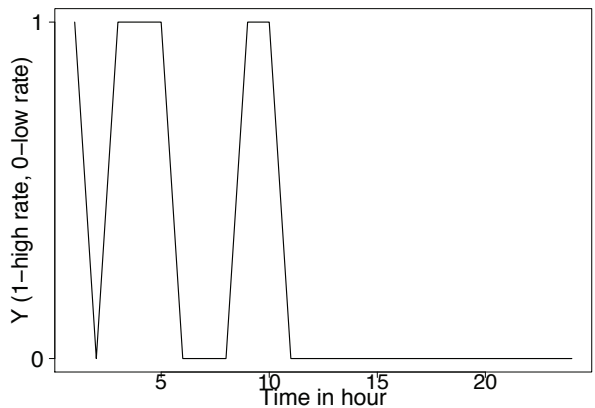}
	\caption{Binary respiratory time series.}
	\label{pre1}
\end{figure}

Of scientific interest is the potential association between self-reported stress and respiratory, or breath, rate measured as the number of breaths per 60 second period.  For the purposes of illustration, we consider a participants breath rate averaged over one-hour intervals starting from midnight and running to midnight over the maximum of a 24  hour period. Empirical data suggests that a respiratory rate of over 20 breaths per minute is considered high for a healthy adult \citep{Barrett:2012}. As such, the time series in this study are discretized into a binary response using a threshold of greater than 20 breaths/min.  Accordingly, if we denote $Y_{it}$ as the average breath rate for subject $i$ at hour $t$, we define $Y_{it}=1$ if the observed average respiratory rate is greater than 20 breaths/min, and 0 otherwise. To illustrate, Fig.~\ref{pre1} presents the observed time series for a randomly sampled participant. Table \ref{transition} depicts the empirical transition table of respiratory rate across all subjects. It illustrates a strong association between the current realization of $Y_{it}$ and lagged values of $Y_{i,t-1}$ and $Y_{i,t-2}.$ In this study, one scientific question of interest is whether or not a potential interaction exists between the lagged realizations $Y_{i,t-1}, Y_{i,t-2}$ and a participant's observed stress level $X_i$. Specifically, it is hypothesized that the association between lagged responses and current breath rate is lower among individuals reporting high stress due to the erratic breathing patterns that high stress situations can evoke. As such, we consider a LARX model including the lagged realization, an indicator for high stress ($1_{[X_i>7]}$), and their interaction. In this study, similar to the discussion in \citet{Holmes:1967}, a subject is considered to be in high stress if the scale exceeds $7$. 
%
\begin{table}[h]
	\centering
	\caption{Empirical transition table of respiratory rate across all subjects}
	\begin{tabular}{ccc}
		\hline
		{Lagged respiratory rate} & \multicolumn{2}{c}{Current respiratory rate}\\
		\cline{2-3}
		&$Y_{i,t}=0$&$Y_{i,t}=1$\\
		\hline

		$Y_{i,t-1}=0$&0.865&0.046\\
		$Y_{i,t-1}=1$&0.038&0.051\\
		$Y_{i,t-2}=0$&0.855&0.047\\
		$Y_{i,t-2}=1$&0.038&0.060\\
		
		\hline
	\end{tabular}
	\label{transition}
\end{table}
\subsection{Fitting the LARX model to the respiratory binary time series data}
We consider LARX(1) and LARX(2) models fitted across the 113 subjects with the same parameter. Stress level $1_{[X_i>7]}$ and the interactions between stress level and past values of the binarized respiratory rate $1_{[X_i>7]} \times Y_{it}$, $1_{[X_i>7]}\times Y_{i,t-1}$ were considered to be covariates. With the independence assumption across subjects, we fit a log likelihood function that is the sum of the log likelihood function (\ref{likelax}) for each subject. Table \ref{table3} provides 95\% confidence intervals for the functionals $\prob(Y_{it}=1\mid Y_{i,t-1})$ and $\frac{\prob(Y_{it}=1\mid Y_{i,t-1})}{\prob(Y_{it}=0\mid Y_{i,t-1})}$ after fitting the LARX(1) model. It can be seen that the confidence intervals derived from Ex-FI are consistently shorter than Em-FI. Specifically, when $Y_{i,t-1}=1$ the confidence interval for $\prob(Y_{it}=1\mid Y_{i,t-1})$ derived from Ex-FI excludes 0.5 (odds excludes 1), while the confidence interval resulting from the use of Em-FI includes 0.5 (odds includes 1). Under the LARX(2) model, the pattern is more obvious. From Table \ref{table4}, it can be seen that comparing the confidence interval from Ex-FI to Em-FI, the average length of all the functionals are relatively smaller. In the most extreme case the Ex-FI derived confidence interval for the odds of high respiratory rate among high stress individuals is approximately 30\% shorter (and excluding 1), when compared to the confidence interval derived using Em-FI. Using the Ex-FI approach, the lagged realizations are determined to be significantly associated with respiratory rate: expectant mothers with low stress level tend to have low rate if their previous realizations are low. In contrast, the wider Em-FI intervals do not rule out a odds of 1 associated with high prior respiratory state among high stress mothers.  
\begin{table}[H]
	\caption{The 95\% confidence intervals of functionals $\prob(Y_{it}=1\mid Y_{i,t-1})$ (Prob) and $ \frac{\prob(Y_{it}=1\mid Y_{i,t-1})}{\prob(Y_{it}=0\mid Y_{i,t-1})}$ (Odds) obtained by fitting the LARX(1) model with stress level and interaction between stress level and past values of the binarized respiratory rate.}
	\centering
	\begin{small}
		\begin{tabular}{llccccc}
			\hline
			&& \multicolumn{2}{c}{Low Stress ($1_{[X_i>7]}=0$)} && \multicolumn{2}{c}{High Stress ($1_{[X_i>7]}=1$)}\\
			\cline{3-4} \cline{6-7}
			\multicolumn{2}{c}{Previous State/Method} &  Prob  &  Odds &&  Prob & Odds  \\ 
			\hline
			$Y_{i,t-1}=0$  \\
			& Ex-FI	& (0.042, 0.061)	& (0.044, 0.065) &  	& 	(0.027, 0.085) &	(0.028, 0.093)   \\
			& Em-FI   	&   (0.042, 0.061)      &  (0.044, 0.065)  &   	& (0.027, 0.085) &    (0.028, 0.093)  \\
			$Y_{i,t-1}=1$\\  
			& Ex-FI    	&    $ \bm{ (0.505, 0.630)}$ 	&    $ \bm{(1.021, 1.701)}$	&   	&	(0.373, 0.731) &    (0.594, 2.724)   \\
			& Em-FI      &    	 $ \bm{(0.498, 0.635)}$      &  $  \bm {(0.998, 1.738) }$	&     	& (0.366, 0.737)&   (0.577, 2.802)    \\
			\hline
		\end{tabular}
	\end{small}
	\label{table3}
\end{table} 
\begin{table}[H]
	\caption{The 95\% confidence intervals of functionals $\prob(Y_{it}=1\mid Y_{i,t-1}, Y_{i,t-2})$ (Prob) and $ \frac{\prob(Y_{it}=1\mid Y_{i,t-1}, Y_{i,t-2})}{\prob(Y_{it}=0\mid Y_{i,t-1}, Y_{i,t-2})} $ (Odds) obtained by fitting the LARX(2) model with stress level and interaction between stress level and past values of the binarized respiratory rate.}
	\centering
	\begin{small}
		\begin{tabular}{llccccc}
			\hline
			&& \multicolumn{2}{c}{Low Stress ($1_{[X_i>7]}=0$)} && \multicolumn{2}{c}{High Stress ($1_{[X_i>7]}=1$)}\\
			\cline{3-4} \cline{6-7}
			\multicolumn{2}{c}{Previous State/Method} & Prob   &  Odds &&  Prob & Odds \\ 
			\hline
			$Y_{i,t-2}=0,Y_{i,t-1}=0$  \\
			& Ex-FI	& (0.044, 0.064)	& (0.046, 0.068) &  	& 	(0.023, 0.081) &	(0.023, 0.088)   \\
			& Em-FI   	&   (0.044, 0.064)      &  (0.046, 0.068)  &   	& (0.023, 0.080) &    (0.023, 0.087)  \\
			$Y_{i,t-2}=1,Y_{i,t-1}=0$\\  
			& Ex-FI    	&     (0.394, 0.553) 	&   (0.651, 1.241)	&   	&	(0.349, 0.851) &    (0.537, 5.701)   \\
			& Em-FI      &    	(0.385, 0.563)       &   (0.626, 1.290) 	&     	& (0.349, 0.851)&   (0.537, 5.707)    \\
			
			$Y_{i,t-2}=0,Y_{i,t-1}=1$  \\
			& Ex-FI	& (0.100, 0.201)	& (0.117, 0.251) &  	& 	(0.017, 0.230) &	(0.017, 0.299)   \\
			& Em-FI   	&   (0.100, 0.210)      &  (0.111, 0.265)  &   	& (0.015, 0.250) &    (0.016, 0.332)  \\
			$Y_{i,t-2}=1,Y_{i,t-1}=1$\\  
			& Ex-FI    	&     (0.670, 0.787) 	&   (2.033, 3.696)	&   	&	$\bm {(0.535, 0.869)}$ &   $\bm{ (1.151, 6.630)}$   \\
			& Em-FI      &    	(0.653, 0.800)       &   (1.878, 4.001) 	&     	& $\bm {(0.469, 0.896)}$&  $\bm{ (0.885, 8.621)}$    \\
			\hline
		\end{tabular}
	\end{small}
	\label{table4}
\end{table}

\section{Conclusion}

We have demonstrated that applying the Em-FI matrix to serially-correlated data 
may lead to undesirable consequences in inference. Such consequences include 
wider confidence intervals (on the average) and thus potentially misleading inferential results. To overcome these limitations, we derived the exact form and an iterative computation formula of the conditional Fisher information matrix for the general logistic
autoregressive model with (without) exogenous covariates (LAR($p$)/LARX($p$)). Although a normality assumption is necessary when the sample size is not large, simulation studies based on the LAR($p$)/LARX($p$) model  demonstrate the advantages 
of Ex-FI over Em-FI in terms of small sample stability, leading to narrower confidence intervals,  on the average, while maintaining false positive rates at or below nominal levels. Numerically, we established the 
convergence of the exact conditional Fisher information and studied the asymptotic behavior as $T$ grows large. Consequently, analysis of the respiratory binary time series data suggests that using Ex-FI may result in greater statistical power when making inference.
In summary, the Ex-FI matrix is recommended over the Em-FI as it provides  greater stability for small time series and equivalent large sample inference. While the derivation of the Ex-FI is non-trivial, it is computationally tractable because it can be obtained iteratively.  The result is a stable estimator that is easily implementable and more stable, particularly for sample sizes less than 200. 

While the proposed approach is promising, there are still potential directions that can be pursued. For instance, the current framework is based on the normality assumption even though the sample size is not too large. As future work, theoretic results on finite sample distribution of maximum likelihood estimates could be established. Moreover, selection of the order $p$ needs to be taken into serious consideration. Motivated by the works of \citet{Kedem:2002} and \citet{Katz:1981}, we may select the optimal lag order $p$ using either 
the Akaike Information Criterion (AIC) or the Bayesian Information Criterion (BIC) 
which are defined to be $AIC(p) = -2\ell(\bm{\hat{\alpha}, \hat{\beta}}\mid \bm{Y}) + 2p$ and $BIC(p) = -2\ell(\bm{\hat{\alpha}, \hat{\beta}}\mid \bm{Y}) + p\log T$ respectively, where $(\bm{\hat{\alpha}, \hat{\beta}})$ is the maximum likelihood estimator of $(\bm{\alpha, \beta}).$





\end{document}